\documentclass{amsart}       
\usepackage[T1]{fontenc}
\usepackage[applemac]{inputenc}
\usepackage[english]{babel}
\usepackage{amssymb}
\newtheorem{theorem}{Theorem}[section]

\newtheorem{proposition}[theorem]{Proposition}

\theoremstyle{definition}

\newtheorem{remark}[theorem]{Remark}
\newtheorem{definition}[theorem]{Definition}

\newtheorem*{convention}{Convention}

\usepackage[isbn=false,doi=false,eprint=false,url=false,style=authoryear,dashed=false]{biblatex}

\makeatletter

\newrobustcmd*{\parentexttrack}[1]{%
  \begingroup
  \blx@blxinit
  \blx@setsfcodes
  \blx@bibopenparen#1\blx@bibcloseparen
  \endgroup}

\AtEveryCite{%
  }

\makeatother

\renewcommand{\cite}{\parencite}
\DeclareNameAlias{sortname}{last-first}
\renewbibmacro*{volume+number+eid}{%
  \printfield{volume}%
  \setunit*{\addnbspace}
  \printfield{number}%
  \setunit{\addcomma\space}%
  \printfield{eid}}
 \DeclareFieldFormat[article]{volume}{\textbf{#1}}
\DeclareFieldFormat[article]{number}{{\mkbibparens{#1}},}
\DeclareFieldFormat{pages}{#1 pages}
\DeclareFieldFormat
  [article,inbook,incollection,inproceedings,patent,thesis,unpublished]
  {pages}{#1}
\renewbibmacro{in:}{}
\newcommand{\ch}{\operatorname{ch}}
\newcommand{\Hom}{\operatorname{Hom}}
\newcommand{\M}{\mathcal{M}}
\newcommand{\MM}{\overline{\mathcal{M}}}
\newcommand{\Dih}{\mathsf{Dih}}   
\newcommand{\Ind}{\mathrm{Ind}}

\newcommand{\sym}{\mathbb{S}}
\newcommand{\F}{\mathbb{M}}
\newcommand{\V}{\mathcal{V}}
\newcommand{\W}{\mathcal{W}}
\newcommand{\E}{\mathcal{E}}
\newcommand{\pp}{\mathfrak{p}}
\newcommand{\qq}{\mathfrak{q}}
\newcommand{\Ass}{\mathsf{Ass}}
\newcommand{\wreath}{\wr}            
\newcommand{\Z}{\mathbb{Z}}

\title{A remark on Getzler's semi-classical approximation}
\author{Dan Petersen}

\begin{document} 

\maketitle

\begin{abstract} Ezra Getzler notes in the proof of the main theorem of ``The semiclassical approximation for modular operads'' that ''A proof of the theorem could no doubt be given using [a combinatorial interpretation in terms of a sum over necklaces]; however, we prefer to derive it directly from Theorem 2.2''. In this note we give such a direct combinatorial proof using wreath product symmetric functions.  \end{abstract}



\section{Introduction}

Let $\V$ be a stable $\sym$-module, i.e.\ a collection $\V((g,n))$ of representations of $\sym_n$ indexed by pairs $(g,n)$ with $2g-2+n > 0$. The paper \cite{getzlerkapranov} defines an endofunctor $\F$ on the category of stable $\sym$-modules, modeled on the way that the moduli spaces $\M_{g,n}$ are glued together to form boundary strata of $\MM_{g',n'}$. In \cite{semiclassical} an explicit formula is derived which describes the genus one part of $\F \V$ in terms of $\V$. The case of genus zero had been described already in \cite{getzler94}, in terms of the Legendre transform. 

If $V$ is a representation of $\sym_n$, let $\ch V$ denote the corresponding symmetric function. The main theorem of \cite{semiclassical} reads\footnote{The term $\frac{1}{4}\psi_2(\mathbf a_0'')$ is missing in Getzler's paper; it was pointed out in \cite{faberconsani} that there is a minor computational error there.}
\[ \mathbf b_1 = \left( \mathbf a_1 - \frac{1}{2}\sum_{n\geq 1} \frac{\phi(n)}{n}\log(1-\psi_n(\mathbf a_0'')) + \frac{\dot{\mathbf a}_0(1+\dot{\mathbf a}_0) + \frac{1}{4}\psi_2(\mathbf a_0'')}{1-\psi_2(\mathbf a_0'')} \right) \circ (h_1 + \mathbf b_0'), \]
where
\[ \mathbf a_g = \sum_{n} \ch \V((g,n)) \]
and
\[\mathbf b_g = \sum_{n} \ch \F\V((g,n))\]
are generating functions. 
 For a symmetric function $f$, $f'$ denotes $\frac{\partial f}{\partial p_1}$ and $\dot f$ denotes $\frac{\partial f}{\partial p_2}$. The $\psi_k$ are the Adams operations defined by $\psi_k(f) = p_k \circ f$. 

Let us introduce some terminology.

\begin{definition}A \emph{graph} is a finite set with a partition and an involution, as in e.g.\ \cite{getzlerkapranov}. A \emph{corolla} is a graph with one vertex. A \emph{necklace} is a graph $\Gamma$ such that $b_1(|\Gamma|) = 1$ and which is not disconnected by removing any edge. \end{definition}

The term $(h_1 + \mathbf b_0')$ can be interpreted combinatorially as taking one copy of the trivial representation, together with all possible graphs corresponding to a stable tree of genus zero vertices with a single distinguished leg. The plethysm should be interpreted as a ``gluing'' operation. In the larger expression 
\[ \left( \mathbf a_1 - \frac{1}{2} \sum_{n\geq 1} \frac{\phi(n)}{n}\log(1-\psi_n(\mathbf a_0'')) + \frac{\dot{\mathbf a}_0(1+\dot{\mathbf a}_0) + \frac{1}{4}\psi_2(\mathbf a_0'')}{1-\psi_2(\mathbf a_0'')} \right),\]
the first term describes corollas of genus one, and the claim is that the rest  is the sum over all possible graphs that are given by a necklace of genus zero vertices. Then the plethysm with $(h_1 + \mathbf b_0')$ gives us the sum over all graphs obtained by attaching genus zero trees (possibly empty, corresponding to the trivial representation) to either the genus one vertex or a necklace, which produces a sum over all stable graphs of genus one, and we recover the definition of $\F$.

Hence the meat of the theorem lies in showing that 
\[-\frac{1}{2}\sum_{n\geq 1} \frac{\phi(n)}{n}\log(1-\psi_n(\mathbf a_0'')) + \frac{\dot{\mathbf a}_0(1+\dot{\mathbf a}_0) + \frac{1}{4}\psi_2(\mathbf a_0'')}{1-\psi_2(\mathbf a_0'')}\]
gives exactly the sum over necklaces of genus zero vertices. It is pointed out in the proof of the theorem that there probably exists a direct combinatorial proof of this fact. However, Getzler deduces it by somewhat involved computations using the more general Getzler-Kapranov formula of \cite{getzlerkapranov} which relates $\V$ and $\F\V$ for all $g$ and $n$, and an explicit representation of the so-called plethystic Laplacian in terms of a formal heat kernel over $\Lambda(\!(\hbar)\!)$.  

In this note we give a combinatorial proof of the fact that the sum over necklaces gives exactly this expression, using only standard facts about wreath product symmetric functions. In particular we are able to give a combinatorial interpretation to the terms in the sum: the first is a sum over all rotational symmetries of the necklaces, and the second is a sum over all symmetries under reflection.

\begin{convention}We consider throughout $\sym$-modules and representations in some fixed symmetric monoidal category $\E$ with finite colimits, additive over a field of characteristic zero. The final assumption allows us to identify $\sym$-modules in $\E$ with symmetric functions, i.e. $K_0([\sym,\E]) \cong K_0(\E) \widehat\otimes \Lambda$,
where $\Lambda$ is the ring of symmetric functions graded by degree, and $\widehat \otimes$ is the completed tensor product. We tacitly omit $\E$ from the notation.
 \end{convention}

\section{Cyclically ordered necklaces}

We start by considering the easier case of necklaces which are equipped with a cyclic ordering. This case is used in the article \cite{cuspformmotives}, and it will serve as motivation for the proof in the unordered case.

\begin{definition} An $\sym$-module $\mathcal V$ is the data of a representation ${\mathcal V}(n)$ of $\sym_n$ for each positive integer $n$. (Usually one would include $n=0$, but it will be slightly more convenient for us not to do so.)\end{definition}

\begin{definition} Let $\V$ and $\W$ be $\sym$-modules. We define their direct sum $\V \oplus \W$ componentwise and their tensor product by
\[ (\V \otimes \W)(n) = \bigoplus_{k+l =n} \mathrm{Ind}_{\sym_k \times \sym_l}^{\sym_n} \V(k) \otimes \W(l).\] This makes the category of $\sym$-modules a symmetric monoidal category.  \end{definition}

\begin{definition} Let $\V$ and $\W$ be $\sym$-modules. The plethysm $\V \circ \W$ is defined by
\begin{equation} \label{plethysm} (\V \circ \W)(n) = \bigoplus_{k=1}^\infty \V(k) \otimes_{\sym_k} (\W^{\otimes k}) (n) \end{equation}
where $(\W^{\otimes k}) (n)$ is considered as an $\sym_k$-module by permuting the factors, i.e.\ via the symmetric monoidal structure on $\sym$-modules. \end{definition}

Let $\Ass$ denote the $\sym$-module defined by 
\[ \Ass(n) =  \mathrm{Ind}_{\Z/n\Z}^{\sym_n}\, \mathbf{1}, \]
where $\mathbf 1$ is the trivial representation, i.e.\ the monoidal unit. Pictorially we think of $\Ass(n)$ as describing corollas with $n$ cyclically ordered input legs, or equivalently, with an embedding in the plane. 

\begin{proposition} \label{earlierprop}The plethysm $\Ass \circ \mathbf a_0''$ is the $\sym$-module describing cyclically ordered necklaces of genus zero vertices. \end{proposition}

\begin{proof}Informally, we think of $\mathbf a_0''$ as corollas of genus zero with two marked legs. We think of the first as the ``clockwise'' one and the second as the ``counterclockwise'' one. There is an evident combinatorial bijection between cyclic necklaces of genus zero vertices and collections of genus zero vertices attached along two marked legs to a corolla with cyclically ordered inputs.

More formally, one can check from the definition of plethysm that one gets the correct result, using that 
\[ \mathbf a_0'' = \sum_{n \geq 3} \ch \mathrm{Res}_{\sym_{n-2}}^{\sym_{n}} \V((0,n)) \]
and that tensoring with $\mathrm{Ind}_{\Z/n\Z}^{\sym_n}\, \mathbf{1}$ is the same as taking coinvariants under the action of $\Z/n\Z$. \end{proof}

Let $\Psi \colon \sym_n \to \Lambda$ be the \emph{cycle map}, defined as
\[ \Psi(x) = \prod_{\sigma \text{ a cycle in $x$}} p_{|\sigma|}. \]
Recall that $\Psi$ induces an isomorphism $\ch \colon R(\sym_n) \to \Lambda^n$ via 
\[V \mapsto \frac 1 {n!}\sum_{x \in \sym_n} \mathrm{Tr}(x\mid V) \Psi(x). \]

\begin{proposition}\label{ind}If $H$ is a subgroup of $\sym_n$, then 
\[ \ch \mathrm{Ind}_H^{\sym_n} \mathbf 1 = \frac 1 H \sum_{h\in H} \Psi(h). \]
\end{proposition}

\begin{proof} See \cite[Chapter 1, Section 7, Example 4]{macdonald}. \end{proof}

\begin{proposition}\label{bb}There is an equality of generating series
\[\sum_{n=1}^\infty \ch \Ass(n) = \sum_{n=1}^\infty \frac 1 n \sum_{d \mid n}\phi(d) p_d^{n/d} = -\sum_{n=1}^\infty \frac{\phi(n)}{n}\log(1-p_n).\] \end{proposition}

\begin{proof}The first equality follows from the preceding proposition, and the second by Taylor expanding and equating coefficients. See also \cite[Example 7.6.2]{getzlerkapranov}. \end{proof}

\begin{proposition}\label{cuspformthm}The sum over all cyclically ordered necklaces is given by 
\[- \sum_{n\geq 1} \frac{\phi(n)}{n}\log(1-\psi_n(\mathbf a_0'')). \]
\end{proposition}

\begin{proof}This follows now by putting together Propositions \ref{earlierprop} and \ref{bb}.\end{proof} 

This is the formula needed in \cite{cuspformmotives}. 

\section{Necklaces and wreath products}

A natural way to compute the sum over necklaces in a combinatorial fashion would be to interpret it, too, as a plethysm.  One might let $\Dih$ denote the $\sym$-module whose $n$th component is spanned by necklaces with $n$ vertices considered up to dihedral symmetry, i.e.\ the $\sym_n$-module $\Ind_{D_n}^{\sym_n} \mathbf 1$, and then consider the plethysm $\Dih \circ \mathbf a_0''$.

This will however not give the right answer, and the basic problem with such an approach is that the action of the dihedral group on the dual graph of a necklace does not factor through the map $D_n \to \sym_n$; indeed, $\sym_n$ just acts by permuting the vertices, but the reflections in $D_n$ should act also by switching which of the two marked legs on each vertex should be ``clockwise'' and ``counterclockwise''.

To incorporate the possibility of having automorphisms which switch the two legs, we will have to work instead with the restriction
\[ \sum_{n \geq 3} \mathrm{Res}_{\sym_2 \times \sym_{n-2}}^{\sym_{n}} \V((0,n)) \]
and consider $D_n$ not as subgroup of $\sym_n$ but of the hyperoctahedral group $\sym_2 \wreath \sym_n = (\sym_2)^n \rtimes \sym_n.$ Let $G$ be a finite group.

\begin{definition} A $(G\times \sym)$-module $\V$ is a sequence $\V(n)$ of representations of $G \times \sym_n$. \end{definition}

\begin{definition} A $(G\wreath\sym)$-module $\W$ is a sequence $\W(n)$ of representations of $G \wreath \sym_n$. \end{definition}

Sums and tensor products of $(G \times \sym)$- and $(G \wreath \sym)$-modules are defined in the same way as for $\sym$-modules. 

\begin{definition} Let $\V$ be a $(G \times \sym)$-module and $\W$ a $(G\wreath\sym)$-module. We define the plethysm $\W \circ_G \V$ by
\[ (\W \circ_G \V)(n) = \bigoplus_{k\geq 0} \W(k) \otimes_{G \wreath \sym_k}{ (\V^{\otimes k})(n)}. \]
Note that if $G$ acts on an object $V$ of a symmetric monoidal category, then $G \wreath \sym_k$ acts on $V^{\otimes k}$, so the tensor product above makes sense. \end{definition}

\begin{remark}\label{deg1}When $\W$ is concentrated in degree 1, then $\W$ is just a representation of $G$ and we recover the ordinary tensor product of $G$-representations, i.e.\
 $\W \circ_G \V = \W \otimes_G \V$.  \end{remark}

We consider the dihedral group $D_n$ as the subgroup of $\sym_2 \wreath \sym_n$ generated by  the elements
$ (\mathbf 1,\tau)$  and  $(-\mathbf 1,\sigma)$
where  $\mathbf 1 \in \sym_2^n$ is the element $(1,1,...,1)$, $-\mathbf 1$ is the element $(-1,-1,...,-1)$, $\tau \in \sym_n$ is the $n$-cycle $(12\cdots n)$, and $\sigma$ is the reflection $(1n)(2,n-1)\cdots$. 

\begin{definition}Let $\Dih$ be the $(\sym_2\wr\sym)$-module defined by $\Dih(n) = \mathrm{Ind}_{D_n}^{\sym_2\wr\sym_n}\mathbf 1$. \end{definition}

\begin{definition}For an $\sym$-module $\V$, let $\V^{(n)}$ denote its restriction to an $(\sym_n \times \sym)$-module. 
\end{definition}

\begin{proposition}\label{dih}The underlying $\sym$-module of the $(\sym_2 \times \sym)$-module
\[ \Dih \circ_{\sym_2} \mathbf a_0^{(2)}\]
is the submodule of $\mathbf b_1$ of  unordered necklaces of genus zero vertices. \end{proposition}

\begin{proof} The proof is now the same as the proof of Proposition \ref{earlierprop}.  \end{proof}

To describe the $\sym_2 \wreath \sym$-module $\Dih$, we shall need to work with the ring $\Lambda(G)$ of \emph{wreath product symmetric functions}. This ring is defined in \cite[Chapter I, Appendix B]{macdonald}. The ring $\Lambda(G)$ is generated as an algebra by generalized power sums $p_n(c)$ where $n$ is a positive integer and $c$ is a conjugacy class of $G$. The degree of $p_n(c)$ is $n$.  There is a natural map
\[ \Psi: G \wreath \sym_n \to \Lambda(G)^n\]
generalizing the cycle map $\sym_n \to \Lambda^n$. One computes $\Psi(g_1,...,g_n, x)$ as follows: for each cycle $\sigma$ of $x$, take the product of the corresponding $g_i$; this product lies in a well-defined conjugacy class $c(\sigma)$ of $G$. Then 
\[ \Psi(g_1,...,g_n, x) = \prod_{\sigma \text{ a cycle in $x$}} p_{|\sigma|}({c(\sigma)}).  \]
As before there is an isomorphism onto the degree $n$ part, 
$\ch \colon  R(G \wreath \sym_n) \to \Lambda(G)^n,$ 
defined by
\[V \mapsto \frac 1 {|G|^nn!}\sum_{x \in G \wreath \sym_n} \mathrm{Tr}(x\mid V) \Psi(x). \]
The plethysm of $(G \wr \sym)$-modules and $(G \times \sym)$-modules can now be described equivalently as an action of $\Lambda(G)$ on $R(G) \otimes \Lambda$.

Proposition \ref{ind} holds true for wreath product symmetric functions --- the proof given in Macdonald's book carries over without changes. Hence we have:

\begin{proposition} \label{indd}Let $H$ be a subgroup of $G \wreath \sym_n$. Then 
\[\Ind_H^{G\wreath\sym_n} \mathbf 1 = \frac{1}{|H|} \sum_{h\in H} \Psi(h) \in \Lambda(G).\]\end{proposition}

\begin{proposition}\label{generating}Let $G = \sym_2$, and denote the power sums in $\Lambda(\sym_2)$ corresponding to the identity conjugacy class by $\pp_n$ and the power sums corresponding to the non-identity by $\qq_n$. Then 
\[ \sum_{n\geq 1} \ch \Dih(n)  = -\frac{1}{2}\sum_{n\geq 1} \frac{\phi(n)}{n} \log(1-\pp_n) + \frac{\frac{\qq_1}{2}(1+\frac{\qq_1}{2}) + \frac{1}{4} \pp_2  }{1-\pp_2}.\]\end{proposition}

\begin{proof}From Proposition \ref{indd} and the definition of $\Psi$ one sees that
\[ \ch \Ind_{D_n}^{\sym_2 \wreath \sym_n}\mathbf 1 = \frac{1}{2n} \sum_{d|n} \phi(d) \pp^{n/d}_{d} +\begin{cases} 
\frac{1}{4}(\qq_1^2 \pp_2^{n/2-1} + \pp_2^{n/2+1}) & \text{$n$ even} \\ 
\frac{1}{2} \qq_1 \pp_2^{(n-1)/2} & \text{$n$ odd} \\ 
\end{cases} \]
where the first term is the sum over all rotations in $D_n$ and the second is the sum over all reflections. Comparing this with the result of Taylor expanding the logarithms and the geometric series gives the result.
\end{proof}

\begin{proposition}\label{deg1sym}Let $G = \sym_n$. The isomorphism $\Lambda(\sym_n)^1 \to \Lambda^n$ is given as follows: if the conjugacy class $c$ in $\sym_n$ is given by the cycle type $(\lambda_1,\ldots,\lambda_k) \vdash n$, then $p_1(c) \mapsto p_{\lambda_1} \cdots p_{\lambda_k}$. \end{proposition}

\begin{proof}One needs only to compare the different isomorphisms 
\[ \Lambda^n \leftarrow  R(\sym_n) = R(\sym_n \wr \sym_1) \to \Lambda(\sym_n)^1. \qedhere \]\end{proof}

For a symmetric function $f(p_1,p_2,\ldots) \in \Lambda$, let $D(f) = f(\frac{\partial}{\partial p_1}, 2\frac{\partial}{\partial p_2},\ldots)$. 

\begin{proposition}\label{getzlerlemma}Let $f \in \Lambda^k = \Lambda(\sym_k)^1$ and $g \in \Lambda$. Then
\[ f \circ_{\sym_k} g^{(k)} = D(f)g. \]\end{proposition}
\begin{proof} Suppose $\ch U = f$ and $\ch V = g \in \Lambda^{n+k}$. Then
\[ U \circ_{\sym_k} V^{(k)} = U \otimes_{\sym_k}V^{(k)} =  \Hom_{\sym_k}(U, \mathrm{Res}^{\sym_{n+k}}_{\sym_k \times \sym_n} V) \]
by Remark \ref{deg1} and since all representations of $\sym_k$ are self-dual. The characteristic of the latter is equal to the right hand side by \cite[8.10]{getzlerkapranov}. One easily extends the result to virtual representations and non-homogeneous $g$. \end{proof}

\begin{proposition} \label{whatwhat}One has that 
$\pp_n \circ_{\sym_2} f^{(2)} = \psi_n(f'')$,
and $\qq_n \circ_{\sym_2} f^{(2)} = 2\psi_n({\dot{f}})$.
\end{proposition}

\begin{proof}Suppose first that $n=1$. Then $\pp_1$ and $\qq_1$ in $\Lambda(\sym_2)^1$ correspond to $p_1^2$ and $p_2$ in $\Lambda^2$ by Proposition \ref{deg1sym}, so by Proposition \ref{getzlerlemma} we have 
\[ \pp_1 \circ_{\sym_2} f^{(2)} = D(p_1^2) f = f''\]
and 
\[ \qq_1 \circ_{\sym_2} f^{(2)} = D(p_2) f = 2\dot f.\]
In general one has $\pp_n \circ_{\sym_2} f^{(2)} = p_n \circ \pp_1 \circ_{\sym_2} f^{(2)} = \psi_n(f'')$, and $\qq_n \circ_{\sym_2} f^{(2)} = p_n \circ \qq_1 \circ_{\sym_2} f^{(2)} = 2\psi_n(\dot f)$. The associativity and the $\lambda$-ring structure on $\Lambda(G)$ used here is most easily seen from the interpretation as polynomial functors, cf.\ \cite{polynomialfunctors}. \end{proof}

\begin{theorem} The sum
\[-\frac{1}{2}\sum_{n\geq 1} \frac{\phi(n)}{n}\log(1-\psi_n(\mathbf a_0'')) + \frac{\dot{\mathbf a}_0(1+\dot{\mathbf a}_0) + \frac{1}{4}\psi_2(\mathbf a_0'')}{1-\psi_2(\mathbf a_0'')}\]
computes the characteristic of the submodule of $\F\V$ spanned by necklaces of genus zero vertices. \end{theorem}

\begin{proof} This follows by putting together Propositions \ref{dih}, \ref{generating} and \ref{whatwhat}. \end{proof}

\printbibliography
\end{document}